\documentclass[a4paper]{acmtrans2e}
\usepackage{graphicx}
\usepackage{amssymb}
\usepackage{comment}
\usepackage[english]{babel}

\newtheorem{lemma}{Lemma}
\newtheorem{corollary}{Corollary}
\newtheorem{theorem}{Theorem}

\makeatletter
\def\subsubsection{\@ucheadfalse
  \@startsection{subsubsection}{3}{\z@}{6pt plus
    1pt}{-8pt}{\reset@font\normalsize\sffamily}}
\makeatother
\def\paragraph{\subsubsection*}

\newcommand{\mm}{\mathfrak{m}}
\newcommand{\qq}{\mathfrak{q}}
\newcommand{\pp}{\mathfrak{p}}
\newcommand{\ttt}{\mathfrak{t}}
\newcommand{\ee}{\mathfrak{e}}

\newcommand{\ww}{\mathfrak{s}}
\newcommand{\sss}{\mathfrak{r}}

\firstfoot{}
\runningfoot{}
\author{Guillaume Chapuy, LIX, \'Ecole Polytechnique
  \and Michel Marcus, LaBRI, Universit\'e Bordeaux 1
  \and Gilles Schaeffer, LIX, CNRS and \'Ecole Polytechnique}
\title{A bijection for rooted maps on orientable surfaces}

\begin{abstract}
The enumeration of maps and the study of uniform random maps have been
classical topics of combinatorics and statistical physics ever since
the seminal work of Tutte in the sixties.  Following the bijective
approach initiated by Cori and Vauquelin in the eighties, we describe
a bijection between rooted maps, or rooted bipartite quadrangulations,
on a surface of genus $g$ and some simpler objects that generalize
plane trees.
Thanks to a rerooting argument, our bijection allows to compute the
generating series of rooted maps on a surface of genus $g$ with
respect to the number of edges, and to recover the asymptotic numbers
of such maps.

Our construction allows to keep track in a bipartite quadrangulation
of the distances of all vertices to a random basepoint. This is an
analog for higher genus surfaces of the basic result on which were
built the recent advances in the comprehension of the intrinsec
geometry of large random planar maps, hopefully opening the way to the
study of a model of continuum random surfaces of genus $g$.
\end{abstract}

\category{Combinatorics}{G.2.1}{Enumeration; Graph theory}

\date{{between 14 Oct 1998 and december 2007}}

\begin{document}

\begin{bottomstuff}
 {G.S. acknowledges financial support from the french ministery of
   research via grants ACI-MD Geocomp and ANR Sada.}
\end{bottomstuff}
\maketitle

\section{Introduction}
\label{sec:in}

This article is concerned with the enumerative and probabilistic study
of \emph{maps on orientable surfaces}: rougly speaking a map of genus
$g$ is a proper embedding of a graph in $\mathcal{S}_g$, the compact
oriented surface of genus $g$ without boundary.  Maps on the sphere
$\mathcal{S}_0$, also known as \emph{planar maps}, have been studied
in enumerative combinatorics ever since the seminal papers of William
T. Tutte in the sixties \cite{Tutte:census}, and our primary
motivation is still to explore the unexpected enumerative properties
of these objects: For instance, the number of rooted planar maps
having $n$ edges is given by the simple formula
\begin{equation}\label{planar}
\vec{q}_{0,n}=\frac{2\cdot 3^n(2n)!}{(n+2)!\,n!},
\end{equation}
and the number of rooted maps of genus $g$ having 
$n$ edges satisfies
\begin{equation}\label{genusg}
\vec{q}_{g,n}\sim  c_g\cdot n^{\frac52(g-1)}\cdot12^n,
\end{equation}
when $n$ goes to infinity with $g$ fixed.  The remarkable linear
dependency in $g$ of the polynomial growth exponent of this asymptotic
formula has been observed as a common pattern of many classes of maps
in combinatorics \cite{BeCa,Gao}. Formula~(\ref{planar}) for the
planar case was recovered also in the physics literature by mean of
the so-called genus expansions of matrix integrals \cite{BIPZ}, which
also led in the 90's to independant derivations of the asymptotic
formula~(\ref{genusg}) \cite[and reference therein]{DFGiZJ}.  Until
now however it has remained an open problem to give a combinatorial
explanation of the linear dependency in the genus \cite{Be:overview}.

The aim of this paper is to describe a bijection
(Theorem~\ref{thm:bijection}) which reformulates and extends to higher
genus surfaces an earlier bijection of Cori and Vauquelin \cite{CoVa}
between planar maps and well labelled trees. In the planar case, a
rerooting argument then immediately leads to Formula~\ref{planar}
(Corollary~\ref{cor:planar}). In higher genus, the combination of our
bijection with a decomposition inspired from Wright's works on labeled
graphs with fixed excess \cite{Wr77,Wr78} leads to a simple combinatorial
derivation of Formula~(\ref{genusg}) (Corollary~\ref{cor:genusg}), and
to a new expression of the generating series for maps of genus $g$
with respect to the number of edges
(Theorem~\ref{thm:genusg}). Classical results of singularity analysis
allow then to recover easily the asymptotic formula~(\ref{genusg}), and
provides a new expression to the constant $c_g$ in front of the formula
(Corollary~\ref{cor:genusg}).

One motivation for the study of maps in statistical physics is the
interpretation of random maps as discrete models of surfaces in the
context of two-dimensional euclidean quantum gravity
\cite{AmDuJo97}.  From a probabilistic point of view, it is
natural to wonder about the geometry of these random discrete
surfaces.  In the planar case our bijection was already used in
\cite{ChSc04} to show that the distance between two random vertices in
a uniform random planar quadrangulation with $n$ faces is of order
$n^{1/4}$, and more generally to study the process describing the
number of vertices at distance $k$ of a vertex.  This line of research
has lead to remarkable progress in our understanding of the geometry
of large random planar maps \cite{BDFG} and of their continuum limit
\cite{MaMo,LG,LGPa}. From this perspective we provide here the basic
building block for the study of higher genus random discrete surfaces.

\section{Maps, quadrangulations and $g$-trees}
In this section we make precise our use of the graph terminology.

\subsection{Graph drawings and maps} 
\begin{figure}
\centerline{\includegraphics[scale=.5]{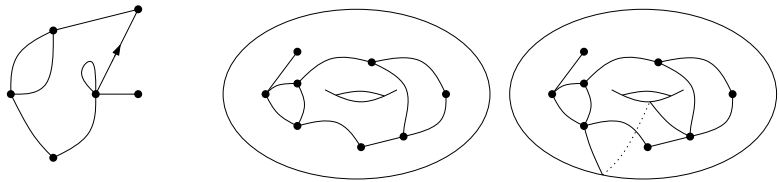}}
\caption{A rooted planar map, a graph drawing on $\mathcal{S}_1$ and a
map of genus $1$. }
\end{figure}
Let $\mathcal S_g$ denote the compact oriented surface of genus $g$
without boundary.  A \emph{graph drawing} on $\mathcal{S}_g$ is a proper
embedding (without edge crossings) of a graph into $\mathcal{S}_g$.
Multiple edges and loops are allowed. A graph drawing on
$\mathcal{S}_g$ is a \emph{map of genus $g$} if moreover each face is
homeomorphic to a disc, or equivalently, if the complement of the
graph in $S_g$ is a union of disjoint simply connected components. In
particular the graph underlying a map is necessarily connected.  Maps
of genus 0 are called \emph{planar maps}.
For enumerative purpose it is convenient to consider \emph{rooted
maps}, that is, pairs $\vec \mm=(\mm,e)$ where $\vec e$ is an arc of
$\mm$ (an arc is an edge with an orientation). The arc $\vec e$ is
called the \emph{root} of $\vec \mm$, its origin is the \emph{root
vertex} of $\mm$ and the face on its right is the \emph{root face}.

In the rest of this text, two maps are considered the same if there
exists an orientation preserving homeomorphism of the embedding
surface that sends one map onto the other (carrying roots if there are
some). These (equivalence classes of) maps can also be defined in a
purely combinatorial way: for instance a rooted map can be described
by listing the local arrangement of arcs around each vertex using any
numbering of these arcs. A modern account of the equivalence between
topological and combinatorial descriptions of maps can be found in
\cite[Ch.~3]{MoTo}.

To any rooted map $\vec\mm$ we associate an arbitrary fixed numbering
of its edges (respectively arcs, vertices, corners, or faces) that
starts with the root edge, which we call the \emph{canonical order} of
the edges (respectively arcs, vertices, corners, or faces) of
$\vec\mm$. This numbering can for instance be given by the first visit
in a depth first search turning counterclockwise around each vertex, or
by any other fixed procedure: in what follows we will only need to
invoke this canonical order to fix notations.

A fundamental result of the theory of surfaces is the Euler
characteristic formula:
\begin{lemma}\label{lem:euler}
If a map of genus $g$ has $v$ vertices, $f$ faces and $n$ edges then
\[
v-n+f=2-2g.
\]
\end{lemma}
The quantity $\chi=2-2g$ is called the Euler characteristic of the map.

\subsection{Maps and bipartite quadrangulations} 
The \emph{degree of a vertex} is the number of edges incident to it,
with the convention that a loop counts for two. The \emph{degree of a
face} is the number of edges incident to it, with the convention that
an edge incident to the same face on both sides counts for two.
\begin{figure}
\centerline{\includegraphics[scale=.5]{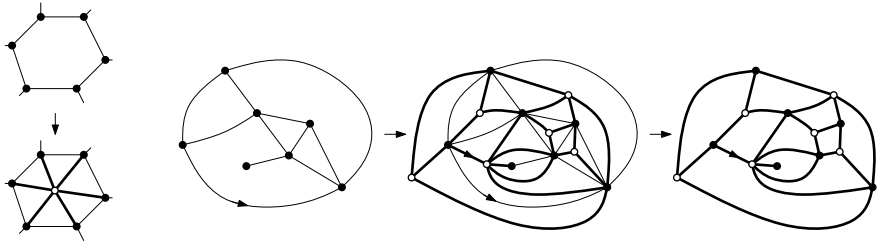}}
\caption{The quadrangulation of a map: local construction and an exemple
in the plane}
\end{figure}
A \emph{quadrangulation} is a planar map having all faces of degree
4. A map is \emph{bipartite} if its vertices are colored in two colors
in such a way that adjacent vertices have different colors (say black
and white). By convention the color of the root vertex of a rooted
bipartite map is always taken to be black.  We recall two standard
results of graph and map theory:
\begin{enumerate}
\item All planar quadrangulations admit a bipartition.
\item For all $n,v,f\geq1$ and $g\geq0$ with $v-n+f=2-2g$, there is a
bijection between maps of genus $g$ with $n$ edges, $v$ vertices and
$f$ faces, and bipartite quadrangulations of genus $g$ with $n$ faces,
$v$ black and $f$ white vertices, and $2n$ edges. Idem between rooted
maps and rooted bipartite quadrangulations.
\end{enumerate}
The first result will not be used but is recalled to stress the fact
that, \emph{a contrario}, a similar statement does not hold for
quadrangulations of genus $g$: there exist non bipartite
quadrangulations of genus $g$ for all $g\geq 1$.

The second result is based on a classical construction which we now
briefly recall: Given a map $\mm$ with black vertices, triangulate
each face from a new white vertex; the new bipartite edges form with
black and white vertices a quadrangulation $\qq$, which is the
bipartite quadrangulation associated to $\mm$. Conversely given a
bipartite quadrangulation $\qq$ create inside each face $f$ a new edge
connecting the two black vertices incident to $f$; these new edges
form with the black vertices a map $\mm$ whose bipartite
quadrangulation is $\qq$.  To get a correspondence in the rooted case,
a rerooting convention must be chosen: if $\mm$ has a root edge $e$,
then let the root of $\qq$ be the arc of $\qq$ with origin the origin
of $e$ and with the face containing $e$ on its right hand side.  This
correspondence is the reason why, although we concentrate in the rest
of this text on bipartite quadrangulations, our results have
implications for all maps.

\subsection{Quadrangulations, distance labelling and $g$-trees}
A \emph{pointed} quadrangulation is a pair $\qq^\bullet=(\qq,v_0)$,
where $v_0$ is a vertex of $\qq$, called the \emph{basepoint}. The
\emph{distance labelling} of $(\qq,v_0)$ is a labelling of the
vertices of $\qq$ by their distance to $v_0$ in the underlying graph:
the label of a vertex $v$ is the minimum number of edges of a path
returning from $v$ to $v_0$.  Observe that in a pointed bipartite
quadrangulation the vertices with the same color as the basepoint are
exactly the vertices at even distance of it. In particular two
adjacent vertices have labels that differ by one in the distance
labelling, and the faces are of two types depending whether the cycle
of incident labels is of the form $(i-1,i,i-1,i)$ or $(i-1,i,i+1)$
for some $i\geq1$.

\begin{figure}
\centerline{\includegraphics[scale=.4]{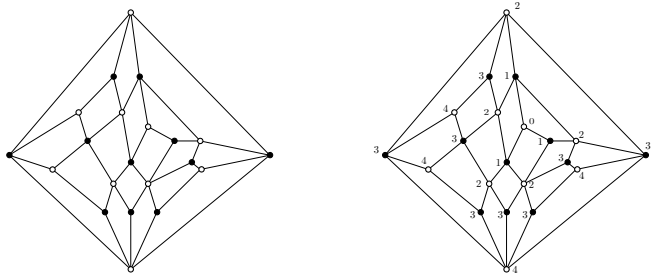}}
\caption{A quadrangulation and its distance labeling}
\end{figure}

A \emph{plane tree} is a planar map with only one face. In particular
a planar tree contains no simple cycle of edges, and the definition of
rooted planar trees agrees with the usual recursive combinatorial
definition of ordered trees: a rooted plane tree can be uniquely
decomposed into a root vertex and a (possibly empty) ordered sequence
of rooted planar trees.  By extension, we call \emph{$g$-tree} a map
of genus $g$ with one face. Observe however that for $g\geq1$, any
$g$-tree contains simple cycles of edges.

\section{Embedded $g$-trees and the opening bijection}
A map is called \emph{embedded} if its vertices are given integer
labels that differ at most by one along every edge, the root vertex
(if the map is rooted) having label 1. A map is \emph{well labeled} if
moreover the support of labels is an interval of the form $(1,m)$. In
other terms a well labeled map is an embedded map with minimum label
$1$.

\subsection{The opening of a pointed bipartite quadrangulation}
\begin{figure}
\centerline{\includegraphics[scale=.48]{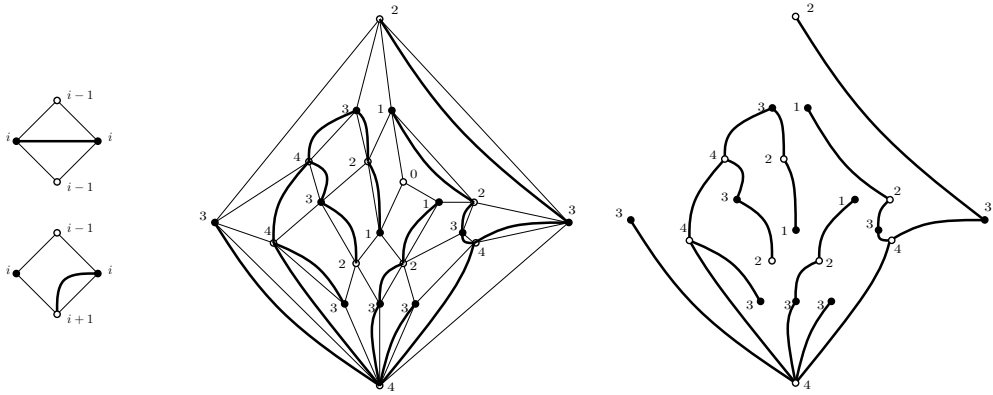}}
\caption{The mapping $\phi$: local rules and an example}\label{fig:rules}
\end{figure}

Given a pointed bipartite quadrangulation $(\qq,v_0)$ of genus $g$
having $n$ faces, we define a mapping $\phi$, the \emph{opening},
which we will later claim produces a well labeled $g$-tree with $n$
edges (see Figure~\ref{fig:rules}):
\begin{itemize}
\item Label vertices of $\qq$ according to their distance to
$v_0$. Recall that the extremities of each edge have labels that
differ by one, so that the labels of vertices on the border of a face
form either a cycle $(i-1,i,i-1,i)$ or a cycle $(i-1,i,i+1,i)$.
\item Create a graph drawing $\phi(\qq,v_0)$ on $\mathcal{S}_g$ with
vertex set the vertices of $\qq$ except $v_0$ and with one edge $e_f$
for each face $f$ of $\qq$, which connects the two corners of $f$ that
have a vertex label larger than their predecessor in clockwise
direction around $f$: in a face with border $(i-1,i,i-1,i)$ the
resulting edge has label $(i,i)$, while in a face with border
$(i-1,i,i+1,i)$ it has label $(i,i+1)$.
\end{itemize}

Obviously the labels of the vertices of $\phi(\qq,v_0)$ are positive
and the variations along its edges belong to $\{+1,0,-1\}$, so that
$\phi(\qq,v_0)$ is well labelled. More interesting is the following
lemma, which proves that $\phi(\qq,v_0)$ is a well labelled $g$-tree.
\begin{lemma}\label{lem:open}
The graph drawing $\phi(\qq,v_0)$ on $\mathcal{S}_g$ has a unique and
simply connected face.
\end{lemma}
\begin{proof} We give here a quick proof which requires some  
familiarity with graph on surfaces, and postpone a more explicit proof
to the appendix. The graph drawing $\qq'=\qq\cup\phi(\qq,v_0)$ is
obviously a map: each face of $\qq$ is divided in $\qq'$ into 2
smaller simply connected faces by an edge of $\phi(\qq,v_0)$ (see
Figure~\ref{fig:rules}). Draw a fake vertex $v_f$ in each face $f$ of
$\qq'$, label it by the smallest label on the border of $f$, and
connect these fake vertices by one fake edge for each edge of $\qq$:
an edge $e$ of $\qq$ with label $(i,i+1)$ that separates faces $f$ (on
its right hand side)

\begin{figure}
\centerline{\includegraphics[scale=.7]{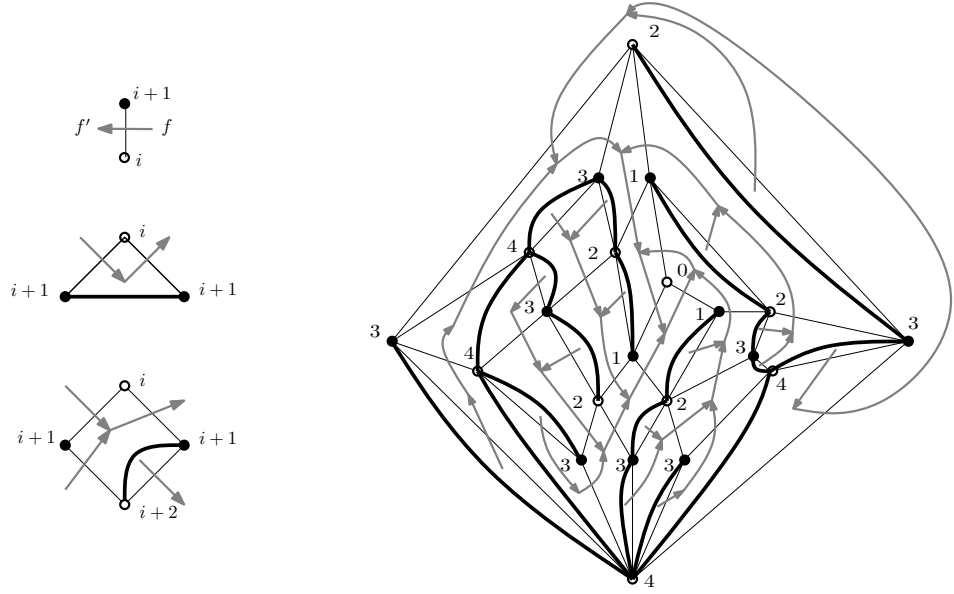}}
\caption{The orientation of dual old edges in $\qq'$: local
configurations and an example.}\label{fig:preuveduale}
\end{figure}

Let us prove that the fake edges form a forest of oriented trees
attached to a unique oriented cycle around $v_0$. Consider the 3 types
of faces of $\qq'$ as represented on Figure~\ref{fig:preuveduale}:
There is exactly one outgoing fake edge in each face, that is, from
each fake vertex, and all oriented fake edges are weakly
decreasing. Hence an oriented cycle of fake edges connects fake
vertices with equal labels. But such a cycle turns counterclockwise
around a unique vertex $v$ of $\qq$, which must be $v_0$ since all its
neighbors have larger labels.

To conclude the proof of the lemma, observe that there can be only one
face in $\phi(\qq,v_0)$ since the fake edges connect all faces of
$\qq'$, and that any loop in this face can be retracted along the tree
like structure of fake edges.
\end{proof}

\subsection{Well labelled polygons, chords diagrams, 
and the closure of a well labelled $g$-tree} We now describe a mapping
$\psi$ which we later claim to be the inverse of $\phi$.  The
construction of $\psi(\ttt)$ will take place in the unique face of
$\ttt$, which is homeomorphic to a disk by definition of maps: in
particular the actual genus of $S_g$ will be somewhat irrelevant for
this reverse construction, which will reveal to be a simple adaptation
of the construction used in \cite{Sc:PhD,ChSc04,MaMo} for the planar
($g=0$) case. More precisely, given a labelled $g$-tree $\ttt$ with
$n$ edges, we decompose it into a rooted $g$-tree $\vec\ttt_0$
(forgetting the labels) and a well labelled rooted $2n$-gon $\vec\pp$
(describing the sequence of labels along the border of the unique face
of $\ttt$).


We will use the following folklore combinatorial result:
\begin{lemma}\label{lem:luka}
Let $\vec\pp$ be a rooted labelled $m$-gon with labels $a_1,\ldots,a_{m}$
in clockwise order starting after the root edge, and assume that the
labels satisfy the \emph{Lukasiewicz} condition $a_{i+1}-a_i\geq -1$,
and $a_i>a_m$ for all $i$.

Let $\pi:\{1,\ldots,m-1\}\to\{2,\ldots,m\}$ be the application that
maps $i$ to its \emph{predecessor}, that is the smallest $j>i$ such
that $a_j=a_i-1$.

Then there exists a unique way to draw all the edges $(i,\pi(i))$
inside $\vec\pp$ so that the resulting graph drawing is a planar map.  In
this map, each inner face is incident exactly once to the polygon,
along an edge $(a_i,a_{i+1})$, and the labels on the border of this
face are then, in clockwise direction,
$(a_i,a_{i+1},a_{i+1}-1,a_{i+1}-2,\ldots,a_{i+1}-\delta)$ where
$\delta=a_{i+1}-a_i+1$ (in particular $\delta\geq0$ and the face incident to
$(a_i,a_{i+1})$ has degree $\delta+2$).
\end{lemma}
\begin{proof} See \cite{Stanley}, or the detailed 
discussion of the appendix.
\end{proof}
\begin{figure}
\centerline{\includegraphics[scale=.5]{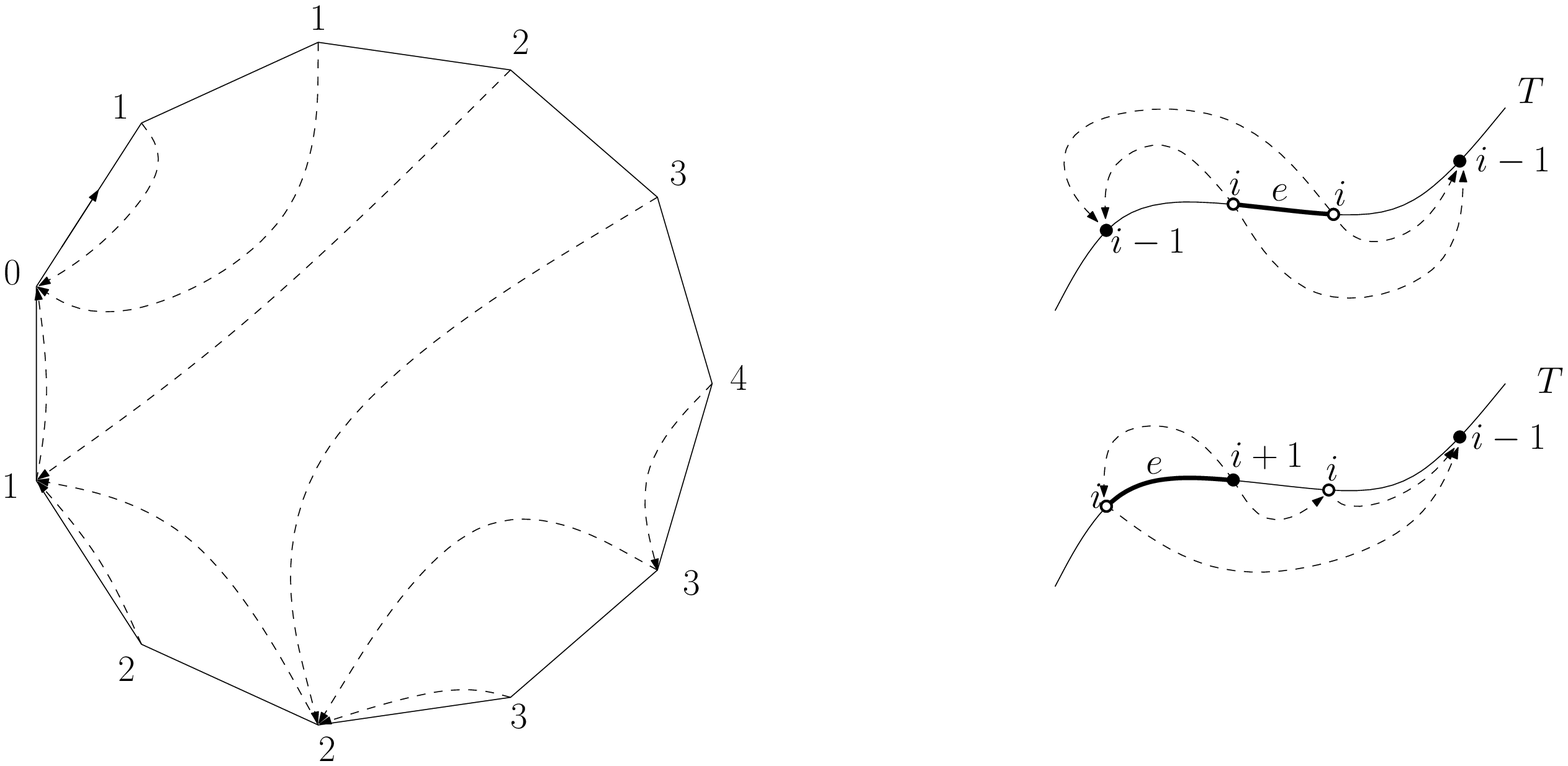}}
\caption{The dissection of a polygon and the two types of faces of
$\psi(\ttt)$}
\label{fig:typesaretes}
\end{figure}

The \emph{closure} of a well labeled $g$-tree $\ttt$ with $n$ edges is
the graph drawing $\psi(\ttt)$ constructed as follows: Let $\pp$ be
the labelled $2n$-gon describing the face $f$ of $\ttt$, and let $k$
be the number of vertices with label $1$ in $\pp$: these vertices
corresponds to \emph{corners} with label $1$ of $f$. Let $\ttt'$ be
the map obtained by adding a vertex $v_0$ with label $0$ in $f$ and
$k$ edges that join this vertex to each of the corners with label $1$
of $f$: the map $\ttt'$ has $k$ faces and the border of each face is a
polygon $\vec\pp_i$ with minimal label~$0$ (taken as root) and labels
that satisfies the Lukasiewicz condition (the variations are even in
$\{+1,0,-1\}$).  In each face of the map $\ttt'$ (on $\mathcal{S}_g$),
that is in each polygon $\vec\pp_i$, the above lemma can thus be used
to add all the edges between corners with positive labels and their
predecessor, in the unique way that provides a map $\ttt''$ on
$\mathcal{S}_g$. The closure $\psi(\ttt)$ of $\ttt$ is the graph
drawing $\ttt''-\ttt$ on $\mathcal{S}_g$ obtained by removing the
edges of $\ttt$ from $\ttt''$.
\begin{lemma}\label{lem:close}
The graph drawing $\psi(\ttt)$ is a quadrangulation of genus $g$ with
$n$ edges.
\end{lemma}
\begin{proof}
According to Lemma~\ref{lem:luka} each face of $\ttt'$ is incident
exactly once to an edge of $\ttt$. Each edge $e$ of $\ttt$ thus
separates two different faces, which are otherwise incident only to
edges of $\ttt'-\ttt$.  Moreover since the variation of labels are
opposite $\pm\delta$ on the two sides of $e$, the removal of $e$
creates a simply connected face of degree $\delta+2+(-\delta+2)-2=4$,
as illustrated by Figure~\ref{fig:typesaretes}. The removal of the $n$
edges of $\ttt$ thus creates $n$ quadrangular faces.
\end{proof}

\subsection{The main theorem}
The main combinatorial result of this paper is the following theorem,
that states that not only $\phi(\qq,v_0)$ is a well labelled map on
$\mathcal{S}_g$ but the mapping $\phi$ is more precisely a bijection.
\begin{theorem}\label{thm:bijection}
The opening $\phi$ is a bijection between 
\begin{itemize}
\item pointed bipartite quadrangulations of genus $g$ that have $n$ faces, $2n$
edges, and $n+2-2g$ vertices ($i$ of which are at odd distance of the
basepoint), and
\item well labeled $g$-trees that have $n$ edges, and $n+1-2g$ vertices
($i$ of them having odd labels),
\end{itemize}
such that the labels of the $g$-tree associated to a bipartite
quadrangulation $\qq$ with basepoint $v_0$ give the distances of
vertices of $\qq$ to $v_0$.
\end{theorem}
\begin{proof}
The fact that $\phi(\psi(\ttt))=\ttt$ follows from the comparison of
Figure~\ref{fig:typesaretes} and Figure~\ref{fig:rules}: for the two
types of faces created in $\psi(\ttt)$, the edge $e$ is indeed
recovered correctly by the opening rules. The fact that
$\psi(\phi(\qq,v_0))=(\qq,v_0)$ follows from the fact that in
$\qq'=\qq\cup\phi(\qq,v_0)$, each edge of $\qq$ joins a corner of
$\phi(\qq,v_0)$ to its predecessor, and from the unicity in
Lemma~\ref{lem:luka}. An alternative direct proof is proposed in the appendix.
\end{proof}

In the enumerative context one usually prefers to deal with rooted
maps to avoid the discussion of symmetries.  The opening is readily
extended to rooted bipartite quadrangulations: define the image of a
rooted quadrangulation $(\qq,\vec e)$ with root vertex $v_0$, to be
$(\phi(\qq,v_0),\vec {e'})$ where $e'$ is the edge created in the root
face of $(\qq,\vec e)$, oriented so that its origin is the endpoint of
$\vec e$ (observe that, according to the rules and the definition of
the root face, the endpoint of $\vec e$ is indeed incident to
$e'$). 
\begin{corollary}
The opening $\phi$ extends to a bijection between
\begin{itemize}
\item rooted bipartite quadrangulations with genus $g$ that have $n$ faces, and
\item rooted well labeled $g$-trees that have $n$
edges.
\end{itemize}
\end{corollary}
\begin{proof}
By construction the origin of $\vec {e'}$ has label $1$ so that
$(\phi(\qq,v_0),\vec {e'})$ is indeed a rooted well labeled tree.
Conversely given a rooted well labeled tree $\vec\ttt=(\ttt,\vec
{e'})$ with root vertex $v$, taking as root edge of $\psi(\vec\ttt)$
the edge $\vec e$ of $\ttt'$ with origin $v_0$ which follows $\vec
{e'}$ in counterclockwise direction around $v$ immediately yields the
unique possible preimage of $\vec\ttt$.
\end{proof}

It turns out that for our purpose an even better variant consists in
taking a basepoint and a root independently. Define the image of a
rooted pointed bipartite quadrangulation $(\qq,\vec e,v_0)$ to be
$(\phi(\qq,v_0),\vec{e'},s)$ where $\vec{e'}$ is the new edge created
in the root face of $\qq$, oriented so that its origin is the
extremity of $\vec e$ with largest label, and where $s\in\{+1,-1\}$ is
a sign recording the orientation of $\vec e$ (say $+$ if $\vec e$ is
increasing with respect to the labeling).  Let moreover $\nu$ be the
relabeling mapping that simultaneously translates all labels of an
embedded map so that the root vertex label is 1. Then


\begin{corollary}\label{cor:embedded}
Opening followed by relabelling is a bijection $\nu\circ\phi$ between
\begin{itemize}
\item rooted pointed bipartite quadrangulations with genus $g$ that have $n$
faces, and
\item pairs consisting of a rooted embedded $g$-tree that
has $n$ edges and a sign $\pm1$.
\end{itemize}
\end{corollary}
\begin{proof}
In view of Theorem~\ref{thm:bijection}, $\phi$ can be used to map
bijectively rooted pointed bipartite quadrangulations $(\qq,v_0,\vec
e)$ of genus $g$ with $n$ faces onto almost well labelled rooted
g-trees $(\phi(\qq,v_0),\vec e)$ with $n$ edges (\emph{almost} well labeled
means the minimum label is 1 but the root vertex has not necessarily
label 1). Now $\nu$ is bijective between almost well labelled
$g$-trees with $n$ edges and embedded rooted $g$-trees with $n$ edges:
indeed given a rooted embedded map (with root vertex label 1 by
definition), there is a unique way to translate its labels so that the
minimum becomes 1.
\end{proof}

The point of this last corollary is that, as we shall see in
Section~\ref{sec:enumeration}, embedded $g$-trees are much easier to
count than well labeled $g$-trees.



\section{The enumeration of quadrangulations 
via schemes}\label{sec:enumeration}

In view of Theorem~\ref{thm:bijection} and
Corollary~\ref{cor:embedded}, the enumeration of maps of genus $g$
according to the number of vertices, faces and edges can be reduced to
the enumeration of rooted well labeled or embedded $g$-trees.  It
turns out that rooted embedded $g$-trees are easier to deal with, so
that we shall rely on the following transposition of
Corollary~\ref{cor:embedded} in terms of generating series:
\begin{corollary}\label{cor:enum}
 The generating series
  \[
  \overrightarrow{Q_g}^\bullet(z)=\sum_{n\geq2g}\vec{q}_{g,n}^{\;\bullet}
  z^n,\quad \textrm{and}\quad T_g(z)=\sum_{n\geq2g}a_{g,n}z^n
  \]
  of rooted pointed bipartite quadrangulations of genus $g$ (with
  respect to the number of faces) and of rooted embedded $g$-trees
  (with respect to the number of edges) are related by the following
  equation:
  \begin{eqnarray*}
    \overrightarrow{Q_g}^\bullet(z)&=&2T_g(z).
  \end{eqnarray*}
\end{corollary}
Since a quadrangulation with genus $g$ and $n$ faces has $n+2-2g$
vertices, $\vec{q}^{\;\bullet}_{g,n}=(n+2-2g)\;\vec{q}_{g,n}$, and the
previous equation can also be rewriten in terms of the generating
series of rooted bipartite quadrangulations of genus $g$ as
\begin{equation}\label{for:enum}
Q_g(z)=\sum_{n\geq2g}\vec{q}_{g,n}z^n=
z^{2g-2}\int \overrightarrow{Q_g}^\bullet(y)\,y^{1-2g}\,dy
\end{equation}
with the initial condition $Q_g(z)=O(z^{2g})$.

\subsection{Planar maps and embedded trees}
Let $T\equiv T(z)=T_0(z)$ be the generating series of rooted embedded
trees with respect to the number of edges. A rooted embedded tree
$\vec\ttt$ which has at least 1 edge can be decomposed into two rooted
trees $\vec\ttt_1$ and $\vec\ttt_2$ by deleting its root edge $e$. The
subtree $\vec\ttt_1$ rooted at the origin of $e$ has root label 1, so
that it is again a rooted embedded tree. The subtree $\vec\ttt_2$
rooted at the endpoint of $e$ has root label $\delta\in\{+1, 0, -1\}$
but, up to a translation of labels by $1-\delta$, it is also a rooted
embedded tree. This decomposition is a bijection
$\vec\ttt\mapsto(\delta,\vec\ttt_1,\vec\ttt_2)$ between non empty
rooted embedded trees with $n$ edges and triples consisting of an edge
variation and two rooted embedded trees with a total of $n-1$
edges. Hence the equation
\begin{eqnarray*}
T(z)=\sum_{\vec\ttt}z^{|\vec\ttt|}=1+3\sum_{(\vec\ttt_1,\vec\ttt_2)}z^{1+|\vec\ttt_1|+|\vec\ttt_2|}=1+3zT(z)^2.
\end{eqnarray*}
The unique formal power series root of this quadratic equation
$T=1+3zT^2$ is
\[
T(z)=\frac{1-\sqrt{1-12z}}{6z}.
\]
In particular $a_n=[z^n]T(z)=\frac{3^n}{n+1}{2n\choose n}$ is the
number of embedded trees with $n$ edges.  For future reference, let us
also observe the relations
\begin{equation}
z\;=\;\frac{T-1}{3T^2},\qquad \qquad
zT'\;=\;\frac{T^2-T}{2-T},
\qquad\textrm{and}\qquad
1+\frac{2zT'}T\;=\;\frac T{2-T}.
\end{equation}
The enumeration of planar quadrangulations then directly follows from
Corollary~\ref{cor:enum}.
\begin{corollary}\label{cor:planar}
The generating series of planar rooted pointed quadrangulations satisfies
\[
    \overrightarrow{Q_0}^\bullet(z)\;=\;2T(z),
\]
and the number of planar rooted quadrangulations with $n$ faces is
\[
\vec q_{0,n}=\frac{\vec{q}^{\;\bullet}_{0,n}}{n+2}=\frac2{n+2}\cdot
\frac{3^n}{n+1}{2n\choose n}.
\]
\end{corollary}

\subsection{Reduced $g$-trees and tree grafting}

We say that a $g$-tree is {\it reduced} if it has no vertices of degree
one.  Upon deleting vertices of degree 1 iteratively, a reduced $g$-tree can
be extracted from any $g$-tree.

\begin{lemma}\label{lem:reduced}
There is a bijection 
between rooted embedded $g$-trees with $n$ edges and pairs consisting
of a rooted embedded reduced $g$-tree $\vec\sss$ with $k$ edges and a forest 
$(\vec\ttt_1,\ldots, \vec\ttt_{2k})$, where,
\begin{itemize}
 \item $\vec\ttt_2,\dots,\vec\ttt_{2k}$ are rooted embedded planar trees.
 \item $\vec\ttt_1$ is either a rooted or a doubly rooted embedded planar tree,
\end{itemize}
and the total number of edges in
$(\vec\sss;\vec\ttt_1,\ldots,\vec\ttt_{2k})$ is $n$.
\end{lemma}
\begin{figure}
\centerline{\includegraphics[scale=0.6]{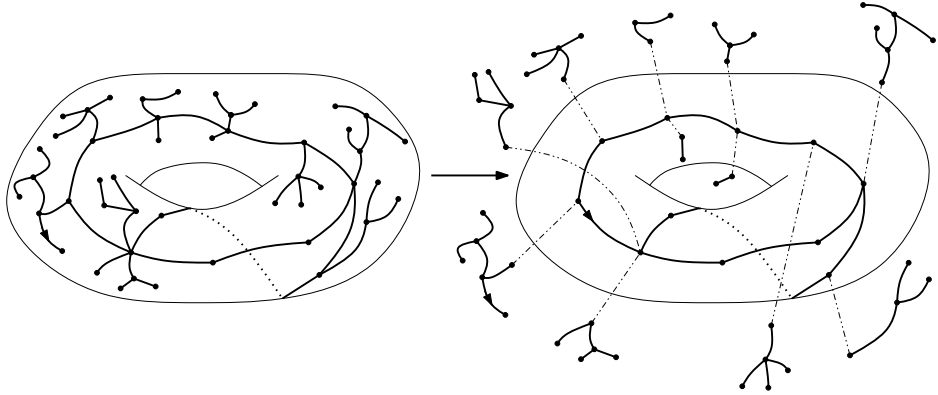}}
\caption{A $1$-tree and its decomposition into a reduced $1$-tree and a
collection of rooted planar trees attached on its corners (only planar
trees which are not reduced to a single vertex are represented).}
\label{fig:reduced}
\end{figure}
\begin{proof}
Given a $g$-tree $\ttt$, a unique reduced $g$-tree $\sss$ is obtained
by iteratively deleting in $\ttt$ the vertices of degree $1$ until
none is left. At the end of this process, the deleted edges form a
forest of rooted embedded planar trees: one (possibly trivial) rooted
tree is attached to each corner of $\sss$ (see
Figure~\ref{fig:reduced}). The root $\vec e$ of $\vec\ttt$, then
provides a root for $\vec\sss$: if $\vec e$ was not deleted it becomes
the root of $\vec\sss$, otherwise $\vec e$ belongs to one of the
deleted trees $\vec\ttt_1$ and the root of $\vec\sss$ is taken to be the
arc of $\sss$ on the right hand side of which $\vec\ttt_1$ is attached.

Once $\vec\sss$ is rooted, its corners are canonically ordered
starting with the corner on the right hand side of the root, and the
deleted trees can be arranged into an ordered list: the $g$-tree
$\ttt$ can then be recovered from $\vec\sss$ and the ordered list
$(\vec\ttt_1,\ldots,\vec\ttt_{2k})$ of trees. Finally, when the root
of $\vec \ttt$ is not an edge of $\vec\sss$ it must be recorded as a
secondary root of $\vec\ttt_1$.
\end{proof}

Lemma~\ref{lem:reduced} immediately yields the following corollary.
\begin{corollary}\label{cor:reduce}
The generating series of rooted embedded $g$-trees $T_g$ and of
rooted embedded reduced $g$-trees $R_g$ are related by:
\begin{eqnarray*}
T_g(z)&=&R_g(zT(z)^2)\cdot \big(1+\frac{2zT'(z)}{T(z)}\big)\;=\;
R_g(zT^2)\cdot \frac T{2-T}.
\end{eqnarray*}
In terms of the primitive $\widehat R_g(t)$, satisfying 
$\frac{t\partial}{\partial t}\widehat R_g(t)= R_g(t)$, the previous
formula reads:
\begin{eqnarray*}
T_g(z)&=&\frac{z\partial}{\partial z}\left(\widehat R_g(zT^2)\right).
\end{eqnarray*}
\end{corollary}
\begin{proof}
The decomposition of a rooted $g$-tree $\ttt$ in a pair
$(\vec\sss;\vec\ttt_1,\ldots,\vec\ttt_{2k})$ where $k$ is the number
of edges of $k$ is a combinatorial composition operation that
results in the composition of generating series:
\begin{eqnarray*}
T_g(z)&=&\sum_{\vec\ttt}z^{|\vec\ttt|}
\;=\;\sum_{\vec\sss;\vec{\vec\ttt}_1,\vec\ttt_2\ldots,\vec\ttt_{2k}\,\mid\,
k=|\vec\sss|}
z^{|\vec\sss|+|\vec{\vec\ttt}_1|+|\vec\ttt_2|+\ldots+|\vec\ttt_{2k}|}\\
&=&\sum_{\vec\sss}z^{|\vec\sss|}\left(\sum_{\vec\ttt}(1+2|\vec\ttt|)z^{|\vec\ttt|}\right)
\left(\sum_{\vec\ttt} z^{|\vec\ttt|}\right)^{2|\vec\sss|-1}
\;=\;R_g(zT(z)^2)\cdot(1+\frac{2zT'}T),
\end{eqnarray*}
where the notation $\vec{\vec\ttt}_1$ stresses the fact that the first
tree can be doubly rooted: given a tree $\ttt$ there are $2|\ttt|$
ways to introduce a double root and one way not to.  The
multiplicative factor $1+2zT'(z)/T(z)$ can also be understood
directly, saying that the first rooted tree (with generating series
$T(z)$) can be replaced by a doubly rooted tree (with generating
series $2zT'(z)$).

The second expression can be checked by a direct comparison with the
first formula:
\begin{eqnarray*}
\frac{z\partial}{\partial z}\left(\widehat R_g(zT^2)\right)&=&
z(T^2+2zTT')\widehat R_g'(zT^2)\\&=&(zT^2)\widehat R_g'(zT^2)(1+2zT'/T)=
R_g(zT^2)\cdot \big(1+\frac{2zT'}{T}\big).
\end{eqnarray*}
\end{proof}

By definition the series $\widehat R_g$ is a generating series of
rooted maps in which rooted maps with $e$ edges are counted with a
weight $1/e$: it is tempting to think about it as a generating series
of unrooted maps, but then each unrooted map is counted with a weight
proportional to the inverse of the size of its automorphism group
(because maps with non trivial automorphisms have less than $2e$
distinct rootings).  In these terms the second expression of
Corollary~\ref{cor:reduce} could be understood directly as follows:
rooted $g$-trees are first constructed as unrooted by grafting trees
in an unrooted reduced $g$-tree, and then rooted. However a rigourous
proof along these lines requires a careful discussion of symmetries,
which is circumvented by working directly with rooted objects.

\subsection{The standard scheme of a reduced $g$-tree}
A \emph{labeled scheme} is a rooted $g$-tree without vertices of degree 1 and
2 that has integer labels on vertices. A scheme is \emph{standard} if
its labels form an integer interval with minimum 0.  Observe that the
labels of a scheme are not required to vary at most by one along
edges, and that in a standard scheme all the integers of the interval
must appear as labels.

Let $\vec\sss$ be a rooted reduced $g$-tree. By definition
the map $\vec\sss$ has no vertices of degree 1. Its vertices of degree
2 are organized into maximal paths, connected together at vertices of
degree at least 3. Let $\vec\ww_0$ and $\vec \ww$ be the maps obtained
from $\vec\sss$ as follows:
\begin{itemize}
\item To get $\vec\ww_0$ from $\vec \sss$, replace each maximal path
made of vertices of degree 2 by an edge, the path containing the root
of $\vec\sss$ providing a rooted edge for $\vec\ww_0$.
\item To get $\vec\ww$ from $\vec\ww_0$, standardize the labels of the
vertices, so that they form an integer interval with minimum label 0:
if the labels appearing on vertices of degree at least 3 in $\vec\sss$
are $\ell_0<\ell_1<\ldots<\ell_p$ then the corresponding labels in
$\vec\ww$ are $0,1,\ldots,p$ (a given label can appear on several
vertices of $\vec\ww_0$, so that $p$ may be smaller than the number of
vertices of $\vec\ww$.).
\end{itemize}
The maps $\vec\ww_0$ and $\vec\ww$ are labeled schemes, and $\vec\ww$
is standard: it is called \emph{the standard scheme of $\vec\sss$}. By
extension, $\vec\ww$ is called the standard scheme of $\vec\qq$ if
$\vec\sss$ was itself extracted from the $g$-tree $\vec\ttt$
associated to a quadrangulation $\vec\qq$. Finally, given a rooted
standard scheme $\vec\ww$ with labels $\{0,\ldots,p\}$ and $p$
positive integers $\ell_1<\ldots<\ell_p$ let
$\vec\ww(\ell_1,\ldots,\ell_p)$ be the rooted labeled scheme obtained
by replacing labels $i$ by $\ell_i$ in $\vec\ww$.

Observe that on every maximal degree 2 path of $\vec\sss$, the labels
perform a sequence of small variations taken from $\{+1,0,-1\}$.  A
\emph{Motzkin walk} $w$ of length $\ell$ is defined as a finite
sequence $w_1\ldots w_\ell$ with $w_i\in\{+1,0,-1\}$. The quantity
$\sum_{i=0}^\ell w_i$ is called the \emph{increment} of $w$. A Motzkin
walk of increment $0$ is called a \emph{Motzkin bridge}. In these
terms, on every maximal path of $\sss$, the sequence of variations of
labels forms a {Motzkin walk}: hence, if the standard scheme $\vec\ww$
of $\vec\sss$ has $k$ edges, a $k$-uple of Motzkin walks
$(\mm_1,\ldots,\mm_k)$ can be extracted from $\vec\sss$, with $\mm_i$
the Motzkin walk read along the path of $\vec\sss$ corresponding to
$i$th edge of $\vec\ww$ in the canonical order.

\begin{figure}
\centerline{\includegraphics[scale=0.6]{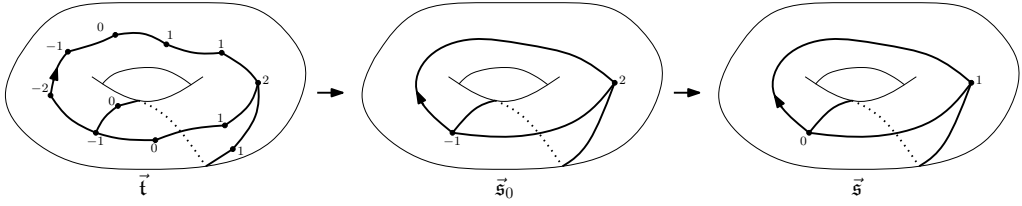}}
\caption{A reduced $1$-tree and its rooted scheme.}
\label{fig:scheme}
\end{figure}

A $k$-uple of Motzkin walks is said to be \emph{compatible} with a
rooted standard scheme with $k$ edges if there exists an embedded
rooted reduced $g$-tree from which they can be both extracted. Observe
that the $g$-tree is unique up to translation of all labels, since it
is recovered by replacing the $i$th edge of the scheme by a chain of
vertices with variations given by the $i$th Motzkin walk.

\begin{lemma}\label{lem:bijection-scheme}
There is a bijection between embedded rooted reduced $g$-trees with
$n$ edges such that the root vertex has degree at least 3, and pairs
formed of a rooted standard scheme of genus $g$ with $k$ edges, and a
compatible $k$-uple of non empty walks with a total of $n$ steps.
\end{lemma}
\begin{proof} 
In view of the definition of compatible $k$-uples of walks, it
suffices to show that if the pair $(\vec\ww;(\mm_1,\ldots,\mm_k))$ is
extracted from embedded rooted reduced $g$-tree $\ttt$ then $\ttt$ can
be recovered from the pair: this is immediate upon replacing the $i$th
edge of $\vec\ww$ by a path with labels given by $\mm_i$, and choosing
the root to be the first edge of the first path.
\end{proof}

In particular, given a rooted standard scheme $\vec\ww$ and a
compatible $k$-uple of Motzkin walks $(\mm_1,\ldots,\mm_k)$ let $\vec
r(\vec\ww;\mm_1,\ldots,\mm_k)$ be the unique corresponding rooted
reduced $g$-tree, and $\vec\ww_0(\vec\ww;\mm_1,\ldots,\mm_k)$ be the
rooted labeled scheme extracted from $\vec
r(\vec\ww;\mm_1,\ldots,\mm_k)$ before standardization.

\begin{lemma}\label{lem:compatible}
Let $\vec \ww$ be a rooted standard scheme with $k$ edges and labels
$\{0,1,\ldots,p\}$, and  let $\ell_1<\ldots<\ell_p$ be positive integers.

Then the set of compatible $k$-uples such that
$\vec\ww_0(\vec\ww;\mm_1,\ldots,\mm_k)=\vec\ww(\ell_1,\ldots,\ell_k)$
is exactly the set of $k$-uples $(\mm_1,\ldots,\mm_k)$ such that, for
all $i$, $\mm_i$ has increment $\ell_{e_i^+}-\ell_{e_i^-}$, where
$e_i^-$ and $e_i^+$ are the labels of the origin and endpoint of the
$i$th edge $\vec e_i$ of $\vec\ww$.
\end{lemma}
\begin{proof}
If a $k$-uple is compatible with $\vec\ww$ and induces
$\ell_0\ldots\ell_p$ then the labels of the corresponding reduced
$g$-tree satisfy the conditions of the lemma. Conversely, given
$\ell_0,\ldots,\ell_p$, a corresponding embedded reduced $g$-tree is
obtained by replacing, for $i=1,\ldots,k$, the $i$th edge $\vec e_i$
of $\vec\ww$ by a labelled path constructed as follows: if
$\mm_i=w_1\ldots w_\ell$ and $e_i^-$ is the label of the origin of
$\vec e_i$ then the inserted path has length $\ell$ and its $j$th
vertex has label $\ell_{e}+\sum_{m=1}^jw_m$, for
$j=0,\ldots,\ell$. The condition on increments ensures that each
vertices of $\vec\ww$ get the same label from all its incident edges.
A root can be chosen arbitrarily on the first path and labels
simultaneously translated to obtain an embedded rooted reduced
$g$-tree. 
\end{proof}

Let $\vec\ww$ be a rooted standard scheme with $v$ vertices, $k$ edges
and labels $\{0,\ldots,p\}$.  We define the \emph{weight} of $\vec\ww$
as the following power series:
\begin{eqnarray*}
W_{\vec\ww}(t) = 
\frac{1}{k}
\sum_{0<\ell_1<\ldots<\ell_{p}}
  \prod_{i=1}^kM_{|\ell_{e_i^+}-\ell_{e_i^-}|}(t),
\end{eqnarray*}
where $e_i^-$ and $e_i^+$ denote the label of the extremities of the
$i$th edge $\vec e_i$ of $\vec \ww$, and, for all $j\geq 0$,
$M_{j}(t)$ is the generating series of non-empty Motzkin walks of
increment $j$ with respect to the length:
$M_j(t)=\sum_{\mm}t^{|\mm|}$, where the sum is over all Motzkin walks
with increment $j$ and length at least 1.

\begin{lemma}
\label{lemma:weight}
The primitive $\widehat{R_g}(t)$ of the generating series of embedded
rooted reduced $g$-trees satisfies
$$
\widehat{R_g}(t) = \sum_{\vec\ww} W_{\vec\ww}(t)
$$ where the sum is taken over all rooted standard schemes of genus $g$. 
\end{lemma}

\begin{proof}
Let us group together in a series $R_{g,k}(t)$ the contributions to the generating
series $R_{g}(t)$ of the embedded reduced $g$-trees whose standard 
scheme has $k$ edges, so that
$R_g(t)=\sum_k R_{g,k}(t)$. Then $2kR_{g,k}(t)$ is the generating series of rooted 
embedded reduced $g$-trees whose standard scheme has $k$ edges and carries a 
secondary root-edge leaving a vertex of degree at least $3$: indeed, each maximal 
path of a given reduced $g$-tree provides $2$ choices of such a secondary root.

Now, these rooted embedded reduced $g$-trees having a secondary root can also be obtained 
by considering first a rooted reduced $g$-tree with a root vertex of degree at least $3$ 
(and whose standard scheme has $k$ edges), and then choosing a root edge. Hence we have:
\begin{eqnarray}
\label{eq:doublerooting1}
2k R_{g,k}(t) = 2\frac{td}{dt}\sum_{|\vec\mathfrak{s}|=k} \widetilde{R}_{\vec\mathfrak{s}}(t)
\end{eqnarray}
where the sum is taken over rooted standard schemes of genus $g$ with $k$ edges, and where 
$\widetilde{R}_{\vec\mathfrak{s}}(t)$ is the generating series of embedded reduced $g$-trees
of scheme $\vec\mathfrak{s}$ with a root vertex of degree at least $3$ 
(the derivative accounts for the choice of the root edge, and the factor
$2$ for the choice of its orientation).

Now, in view of Lemmas~\ref{lem:bijection-scheme}--\ref{lem:compatible}, for each scheme 
$\vec\mathfrak{s}$ with $k$ edges, the 
generating series $\widetilde R_{\vec\mathfrak{s}}(t)$ 
is obtained by summing over all relabellings and compatible $k$-uples
of Motzkin paths. Therefore:
\begin{eqnarray}
\label{eq:doublerooting2}
\widetilde R_{\vec\mathfrak{s}}(t)= 
\sum_{0<\ell_1<\ldots<\ell_{p}}
  \prod_{i=1}^kM_{|\ell_{e_i^+}-\ell_{ e_i^-}|}(t).
\end{eqnarray}
Putting Equations~\ref{eq:doublerooting1} and ~\ref{eq:doublerooting2} together and using the
definition of $W_{\vec\mathfrak{s}}$ gives:
$$
R_{g,k}(t) = \frac{td}{dt}\sum_{|\vec\mathfrak{s}|=k} W_{\vec\mathfrak{s}}(t).
$$
Summing on $k$ and comparing with the definition of $\widehat{R_g}(t)$ yields the result.
\end{proof}

\subsection{The algebra of Motzkin walks}

A Motzkin walk is \emph{primitive} if it has increment $-1$ but none
of its prefix has increment $-1$. In other terms a Motzkin walk
$w=w_1\ldots w_\ell$ is primitive if it satisfies the \L ukasiewicz
property: $\sum_{i=1}^{\ell}w_i=-1$ but for all $j<\ell$,
$\sum_{i=1}^{\ell}w_i\geq0$.
 
Let $U\equiv U(t)=\sum_{w}t^{|w|}$ be the generating series of
primitive Moztkin walks according to the length, and $B\equiv
B(t)=M_0(t)$ be the generating series of non-empty bridges.  Using a
decomposition at the first step for $U$, and at the first return to
the $x$-axis for $B$, one gets:
\begin{eqnarray*}
U&=&t(1+U+U^2)\\
B&=&t(1+2U)(1+B)\;=\;\frac{t(1+2U)}{1-t(1+2U)}.
\end{eqnarray*}
Then the explicit value of $U$ and $B$ are
\begin{eqnarray*}
U=\frac1{2t}(1-t-((1-2t-3t^2)^{1/2}),&\textrm{ and } &
1+B=\frac1{(1-2t-3t^2)^{1/2}}.
\end{eqnarray*}
In terms of $U$, the series $B$ satisfies
\begin{eqnarray*}
B=\frac{U(1+2U)}{1-U^2},&\textrm{ and }&
1+B=\frac{1+U+U^2}{1-U^2}.
\end{eqnarray*}
Decomposing walks at their last passage time at each positive
integer, one obtains the generating function of Motzkin walks 
of increment $i>0$:
\begin{eqnarray*}
M_i(t) = (1+B) U^i
\end{eqnarray*}

The algebra of Laurent power series generated by the series $U$ contains
$t$ since
\[
t=\frac{U}{1+U+U^2}=\frac1{U^{-1}+1+U}.
\]
More generally we would like to characterize rational expressions in
$U$ that are in fact rational expressions in $t$. Obviously these
expressions are symmetric in the exchange $U\leftrightarrow U^{-1}$.
This is also a sufficient condition since the polynomials
$(U+1+U^{-1})^k$ form a base of the set of polynomials in
$\{U,U^{-1}\}$ that are symmetric in these two indeterminates.

\begin{lemma}
A rational series in $U$ is a rational series in $t$ if and only if it
is symmetric in the exchange $U\leftrightarrow U^{-1}$.
\end{lemma}

\subsection{The torus}
Let $\sss$ be a reduced $1$-tree with $k$ edges, and let $d_i$ be the
number of vertices of degree $i$ in $\sss$. Then the vertex-edge
incidence relation and Euler's formula give the following 2 relations:
\[
\left\{
\begin{array}{rcl}
2k &=&\sum_{i\geq1} id_i,\\
k&=&1+\sum_{i\geq1} d_i.
\end{array}
\right.
\]
Eliminating $k$ and recording that $d_1=0$ since $\sss$ is reduced,
this yields
\[
\sum_{i\geq2} (i-2)d_i=2.
\]
Hence, apart from vertices of degree 2, there can be only one vertex
of degree 4 ($d_3=0$, $d_4=1$) or two vertices of degree 3 ($d_3=2$,
$d_4=0$). Let us compute the standard schemes and weights associated to these
two types of reduced $1$-trees:
\begin{description}
\item[1 vertex of degree 4] If the 1-tree consists of 1 vertex of
degree 4 and vertices of degree 2, its standard scheme has just 1 vertex of
degree 4 and 2 edges that are loops. There is a unique rooted $1$-tree
$\vec\ttt_1$ with one vertex of degree four, and giving label 0 to its
unique vertex we obtain the unique rooted standard scheme $\vec\ww_1$ with 1
vertex.  The associated weight is therefore
\[W_{\vec\ww_{1}}=\frac{1}{2}M_0^2=\frac12B^2
=\frac12\frac{(1+2U)^2U^2}{(1-U)^2(1+U)^2}.\]

\item[2 vertices of degree 3] If instead the 1-tree consists of two
vertices of degree 3 and vertices of degree 2, its standard scheme has two
vertices of degree 3 and 3 edges connecting them. There is again a
unique rooted $1$-tree with 2 vertices of degree 3 but the choice of a
label $a$ for the root vertex and $b$ for the second vertex yields
three different rooted standard schemes $\vec\ww_2,\vec\ww_3,\vec\ww_4$:
\begin{itemize}
\item $a=b=0$: in this case the endpoints of all three edges have
label 0, so that the associated weight is 
\[W_{\vec\ww_{2}}=\frac{1}{3}B^3\;=\;
\frac13\frac{(1+2U)^3U^3}{(1-U)^3(1+U)^3}.
\]
\item $a=0$, $b=1$: in this case the associated weight is
\begin{eqnarray*}
W_{\vec\ww_{3}} &=&
\frac{1}{3} \sum_{\ell_1>0} \left((1+B)U^{\ell_1} \right)^3 
\;=\;\frac13(1+B)^3\frac{U^3}{1-U^3}\;=\;
\frac23\frac{(1+U+U^2)^2U^3}{(1-U)^4(1+U)^3}
\end{eqnarray*}
\item $a=1$, $b=0$: this case yields the same weight as the previous one,
$W_{\vec\ww_{4}}=W_{\vec\ww_3}$.
\end{itemize}
\end{description}
Summing over the different cases, we get
\begin{eqnarray*}
\widehat{R_1}&=& W_{\vec\ww_{1}} +
  W_{\vec\ww_{2}} +W_{\vec\ww_{3}}+W_{\vec\ww_4} \\
&=&
\frac{U^{-1}+4+U}{2((1-U)(1-U^{-1}))^2(1+U)(1+U^{-1})}.
\end{eqnarray*}
We thus obtained $\widehat{R_1}$ as a rational function of $U$. Its
symmetry w.r.t. the exchange $U\leftrightarrow U^{-1}$ shows that it is in
fact a rational function of $t$: indeed one easily checks that
\[
\widehat R_1\;=\;\frac{t^2(1+3t)}{2(1-3t)^2(1+t)}.
\]
In view of Corollary~\ref{cor:reduce} we obtain the generating
function of rooted $1$-trees as
\[
T_1(z)=\frac{z\partial}{\partial z}\widehat R_1(zT(z)^2)
\]
and deduce from Corollary~\ref{cor:enum} an expression for the
generating series of pointed rooted bipartite quadrangulation on the
torus as $\overrightarrow{Q_g}^\bullet(z)=2T_1(z)$.  However, in view
of Formula~\ref{for:enum}, one can derive directly the generating
series of rooted bipartite quadrangulations on the torus from
$\widehat R_1$:
\begin{eqnarray*}
Q_1(z)&=&\int \frac{dy}{y}\, \overrightarrow{Q_1}^\bullet(y)\;=\;
2\int \frac{dy}{y}\frac{yd}{dy}\widehat R_1(yT(y)^2)
\\
&=&2\widehat R_1(zT(z)^2)\;=\;2\widehat R_1({\textstyle\frac{T(z)-1}3})
\;=\;\frac{(T-1)^2T}{3(2-T)^2(2+T)},
\end{eqnarray*}
in agreement with the result of \cite{BeCa}.

\subsection{General surfaces} 
In general the incidence and Euler relations for reduced $g$-trees
with $k$ edges and $d_i$ vertices of degree $i$ are
\[
\left\{
\begin{array}{rcl}
2k &=&\sum_{i\geq1} id_i,\\
k&=&2g-1+\sum_{i\geq1} d_i.
\end{array}
\right.
\]
Eliminating $k$ and using $d_1=0$, this yields
\[
\sum_{i\geq2}(i-2)d_i=4g-2.
\]
There is thus again a finite number of possible combinations for the
$d_i$ with $i\geq3$: in particular the two extremal possibilities in
terms of number of edges are:
\begin{itemize}
\item $d_{4g}=1$, and $d_i=0$ otherwise, which leads to schemes with 1
vertex of degree $4g$ and $2g$ edges, the minimal number of edges for
a scheme of genus $g$.
\item $d_3=4g-2$, and $d_i=0$ otherwise, which leads to schemes with
$4g-2$ vertices of degree 3 and $6g-3$ edges, the maximal number of
edges for a scheme of genus $g$.
\end{itemize}

Let $\vec\ww$ be a rooted standard scheme with $q$ vertices and $k$ edges
($k=2g-1+q$), Let us call $p$ the number of distinct nonzero labels in
$\vec\ww$ (in particular $0\leq p\leq q$), $\ee$ the set of edges with
equal extremities, and $\ee'$ its complement.  Then we have by
definition:
\begin{eqnarray*}
W_{\vec\ww}&=&\frac1k B^{|\ee|}\sum_{0<\ell_1<\ldots<\ell_{p}}
\prod_{e\in \ee'}M_{|\ell_{e_-}-\ell_{e_+}|},
\end{eqnarray*}
where $e_-$ and  $e_+$ denote the vertex labels of the extremities
of the edge $e$ of $\ee'$.
 
Observe that for $0\leq i<i'\leq v_{\neq}$, we have
$\ell_{i}<\ell_{i'}$ and
\[
M_{\ell_{i'}-\ell_{i}}=(1+B)U^{\ell_{i'}-\ell_i}=(1+B)\prod_{j=i+1}^{i'}U^{\ell_{j}-\ell_{j-1}},
\]
so that 
\begin{eqnarray*}
W_{\vec\ww}&=&\frac1k B^{|\ee|}(1+B)^{|\ee'|}
\sum_{0<\ell_1<\ldots<\ell_{p}}
\prod_{e\in \ee'}\prod_{j=e_-+1}^{e_+ }U^{\ell_{j}-\ell_{j-1}}.
\end{eqnarray*}
Since the expression summand only involves the differences
$\delta_j=\ell_{j}-\ell_{j-1}$, the sum can be rewritten in terms of
these $p$ new variables: 
\begin{eqnarray*}
W_{\vec\ww}&=&\frac1k B^{|\ee|}(1+B)^{|\ee'|}
\sum_{\delta_1,\ldots,\delta_{p}}
\prod_{e\in \ee'}\prod_{j=e_-+1}^{e_+}U^{\delta_j},
\end{eqnarray*}
where the $\delta_j$ are independant summation indices over positive
integers.  Now for each couple $(j,e)$ such that ${e_-}< j\leq e{_+}$
a contribution $U^{\delta_{j}}$ is obtained. Let $d(j)$ denote the
number of edges $e$ such that the couple $(j,e)$ satisfies the previous
condition. Observe that the $d(j)$ only depend on $\vec\ww$, not on the
$\delta_j$. Then
\begin{eqnarray*}
W_{\vec\ww}&=&\frac1k B^{|\ee|}(1+B)^{|\ee'|}
\sum_{\delta_1,\ldots,\delta_{v_p}}
\prod_{j=1}^{p}U^{d(j)\delta_j}\\
&=&\frac1kB^{|\ee|}(1+B)^{|\ee'|}
\prod_{j=1}^{p}\frac{U^{d(j)}}{1-U^{d(j)}}.
\end{eqnarray*}
Setting $e_==|\ee|$ and $e_{\neq}=|\ee'|$, and expressing $B$ in terms of
$U$ we obtain
\begin{eqnarray*}
W_{\vec\ww}
&=&\frac1k\frac{U^{e_=}(1+2U)^{e_=}(1+U+U^2)^{e_{\neq}}}{(1-U^2)^{k}}
\prod_{j=1}^{p}\frac{U^{d(j)}}{1-U^{d(j)}}.
\end{eqnarray*}

In view of Lemma \ref{lemma:weight}, the generating series of reduced labeled 
$g$-trees is obtained upon summing these individual contributions for 
all rooted standard schemes of genus $g$.

\begin{theorem}\label{thm:genusg}
Let $\mathcal{W}_g$ denote the (finite) set of rooted standard schemes of genus
$g$. Given such a standard scheme $\vec\ww$ with $k$ edges and labels
$\{0,1,\ldots,p\}$, let $e_=$ denote the number of edges with equal
labels at both ends, and $e_{\neq}=k-e_=$. For $j=1,\ldots,p$, denote
moreover by $d(j)$ the number of edges $e$ with endpoint labels
satisfying $e_-<j\leq e_+$, and $d=\sum_{j=1}^pd(j)$.  Then
\begin{eqnarray}
\label{eq:WgU}
\widehat R_{g}(t)
&=&
\sum_{\vec\ww\in\mathcal{W}_g}
\frac1{k}\frac{U^{d+e_=}(1+2U)^{e_=}(1+U+U^2)^{e_{\neq}}}
{(1-U)^{k+p}(1+U)^{k}} 
\prod_{j=1}^{p}\frac1{1+U+\ldots+U^{d(j)-1}}
\end{eqnarray}
where $U$ is the formal power series in $t$ satisfying $U=t(1+U+U^2)$. Moreover
\[
   T_g(z)=\frac{z\partial}{\partial z}\widehat R_g(zT^2),\qquad
    \overrightarrow{Q_g}^\bullet(z)=2T_g(z)
\]
where $T$ is the formal power series in $z$ satisfying $T=1+3zT^2$.
\end{theorem}
In particular we expressed $\overrightarrow{Q_g}^{\bullet}(z)$ as a
rational function of $U(zT^2)$, which is an algebraic function of
degree 4. It is proved however in \cite{BeCa} that $Q_g(z)$ is in fact
a rational function in $z$ and $\sqrt{1-12z}$.  One can verify that
this result implies that our series $\widehat R_g(t)$ is in
fact a rational function in $t$, or in other terms that the expression
of $\widehat R_g(t)$ in terms of $U$ is in fact symmetric in the
exchange $U\leftrightarrow U^{-1}$. We are currently unable to give a
combinatorial proof of this fact.

\subsection{Asymptotics}

In Equation~(\ref{eq:WgU}), there are a finite number of critical
values of $U$, all being roots of unity, and $U=1$ is one of them.  As
a fuction of $t$, $U$ reaches the value $U=1$ at its singularity
$t=1/3$. Moreover, $U$ being a power series in $t$ with positive
coefficients, one has $\left|U(t)\right|\leq 1$ for every complex
number $t$ with $|t|\leq 1/3$, with equality only at $t=1/3$. This
implies that $t=1/3$ is the only dominant singularity of $R_g(t)$.
Finally, $U(t)$ is an algebraic function, which admits a Puiseux expansion
at its singularity:
\[
U=1-\sqrt{3}\sqrt{1-3t}+O(1-3t).
\]
This implies the following corollary:

\begin{corollary}
The dominant singularity of the algebraic function $\widehat R_g(t)$ is at
$t=1/3$, around which the following Puiseux expansion holds:
\[
\widehat R_g(t)=\sum_{\vec\ww\in\mathcal{W}_g}
\frac1{k}\left(\frac32\right)^{k}
\left(\frac{1}{\sqrt3\sqrt{1-3t}}\right)^{k+p}
\prod_{j=1}^{p}\frac1{d(j)}
\;\cdot\;(1+O(\sqrt{1-3t})),
\]
where for each $\vec\ww$ in the summation, $k$ is the number of
edges and $\{0,1,\ldots,p\}$ the set of distinct labels.
\end{corollary}
The dominant terms in this sum are given by the standard schemes $\vec\ww$ with
$k+p$ maximal: $p$ is at most $q-1$ if $\vec\ww$ has $q$ vertices and
all labels are distinct. Then $k+q-1=2q+2g-2$, which is maximal for
$q=d_3=4g-2$, resulting in $k+q-1=10g-6$. 
\begin{corollary}
Let $\mathcal{W}^m_g$ denote
the set of rooted standard schemes of genus $g$ with $4g-2$ vertices of degree
3 that are labelled by distinct integers. Then
\[
\widehat R_g(t)=
\frac1{6g-3}
\frac{3^g}{2^{6g-3}}
\left(
\sum_{\vec\ww\in\mathcal{W}^m_g}\prod_{i=1}^{4g-3}\frac1{d(i)}
\right)
\left(\frac{1}{{1-3t}}\right)^{5g-3}
\;\cdot\;(1+O(\sqrt{1-3t})).
\]
\end{corollary}

In order to obtain the asymptotic behavior of
$\overrightarrow{Q_g}^\bullet(z)$, we need to consider the behavior of
$T$: this algebraic series has a unique dominant singularity at $z=1/12$ and
 \[T=2-2\sqrt{1-12z}+O(1-12z),\] so that
\[zT^2=\frac{T-1}3=\frac13-\frac2{3}\sqrt{1-12z}+O(1-12z).\]
Since $T$ is a power series in $z$ with positive coefficients, one has
$\left|zT^2\right|<1/3$ for every $z\neq1/12$ such that $|z|\leq
1/12$, so that we can compose the singular expansions of $T$ and
$\widehat R_g$ in the expression
$\overrightarrow{Q_g}^\bullet(z)=\frac{zd}{dz}\widehat R_g(zT^2)$ of
Corollary~\ref{cor:enum} to obtain:
\begin{corollary}
The algebraic series $\overrightarrow{Q_g}^\bullet(z)$ has a unique
dominant singularity at $z=1/12$, around which the following Puiseux
expansion holds:
\begin{eqnarray*}
\label{eq:devQ_g}
\overrightarrow{Q_g}^\bullet(z)&=&\frac{5g-3}{6g-3} \frac{3^g}{2^{11g-6}}
\left(
\sum_{\vec\ww\in\mathcal{W}^m_g}\prod_{i=1}^{4g-3}\frac1{d(i)}
\right)
\left(\frac{1}{\sqrt{1-12z}}\right)^{5g-1}\cdot(1+O((1-12z)^{1/4})\\
\end{eqnarray*}
\end{corollary}

It is a well known fact that algebraic functions are amenable to singularity
analysis (see \cite{FlOd,FlSe}). In other words, classical transfer theorems
apply to the previous expansion, leading to:
\begin{corollary}\label{cor:genusg}
The number $\vec q^{\,\bullet}_{g,n}$ of rooted pointed bipartite
quadrangulations of genus $g$ with $n$ faces satisfies:
\begin{eqnarray*}
\vec
q^{\,\bullet}_{g,n}&\sim&\frac{1}{\Gamma\left(\frac{5g-1}2\right)}
\frac{5g-3}{6g-3}\frac{3^g}{2^{11g-6}}
\left( \sum_{\vec\ww\in\mathcal{W}^m_g}\prod_{i=1}^{4g-3}\frac1{d(i)} \right) \cdot
{n^{\frac{5g-3}{2}}}\cdot 12^n.
\end{eqnarray*}
\end{corollary}

Recalling that $ q^\bullet_{g,n} = (n+2-2g) \vec q^{\,\bullet}_{g,n} $
and using the product relation
$\Gamma\left(\frac{5g-1}2\right)=\frac{5g-3}2\Gamma\left(\frac{5g-3}2\right)$,
we finally have:
\begin{corollary}
The number of rooted bipartite quadrangulations of genus $g$ with $n$ 
faces satisfies
\begin{eqnarray*}
q_{g,n}&\sim&\frac{3^g}{(6g-3)2^{11g-7}\Gamma\left(\frac{5g-3}2\right)}
\left(
\sum_{\vec\ww\in\mathcal{W}^m_g}\prod_{i=1}^{4g-3}\frac1{d(i)}
\right)
\cdot
{n^{\frac52(g-1)}}\cdot 12^n.
\end{eqnarray*}
\end{corollary}

Observe that the maps in $\mathcal{W}^m_g$ are vertex labeled rooted
maps with $4g-2$ vertices of degree 3, one face and degree $g$. Their
number is $(4g-2)!\varepsilon_g$ where $\varepsilon_g$ is the number
of cubic maps of genus $g$ with one face (or dually triangulations
with one vertex), known to be \cite{WaLe}:
\[
\varepsilon_g 
=\frac{2}{12^{g}}\frac{(6g-3)!}{g!(3g-2)!}.
\]
In particular this allows to give bounds on the constant
$\sum_{\vec\ww\in\mathcal{W}^m_g}\prod_i\frac1{d(i)}$. For fixed $g$, this
constant is in principle a computable finite sum.  However, in
practice, since the number of maps in $\mathcal{W}^m_g$ grows
superexponentially in $g$, the computation is difficult to perform in
a reasonable time. In particular we were able to obtain in this way
the explicit value of the constants only for $g=1,2$, leading to the
following asymptotics, in agreement with \cite{BeCa}:
$$q_{1,n} \sim \frac{1}{24}\cdot 12^n \;\mbox{ and }\;
 q_{2,n}\sim\frac{7}{4320\sqrt{\pi}}\cdot n^{5/2}\cdot 12^n.$$ Several
 different explicit non linear recursions are known for these
 constants (see \cite{BeCa,DFGiZJ,Witten}), but we were unable to
 derive them from our expressions.

\section{Conclusion}

As in the planar case, our main bijection maps the distances to the
basepoint of a quadrangulation onto the labels of vertices of the
associated $g$-tree.  In particular it allows to prove that the
distances of a random vertex to the basepoint in a uniform random
quadrangulation with $n$ faces is again of order $n^{1/4}$: a detailed
analysis in this direction was recently done in \cite{Mi07} for
the Bolzmann distribution on quadrangulations of genus $g$. A natural
problem, partially considered there, is more generally to construct a
continuum limit of large quadrangulations in this $n^{1/4}$ scaling:
in view of the results of \cite{LG,LGPa} in the planar case, a
continuum limit of quadrangulations of genus $g$ is expected to have
the topology of $\mathcal S_g$.

We conjecture that the resulting \emph{continuum random
quadrangulation of genus $g$} would be universal enough to describe
the continuum limit of all models of uniform random maps with simple
face degree or connectivity constraints, but also the continuum limit
of uniform random graphs with $n$ vertices and minimum genus $g$.
\medskip

As we have seen the dominant terms in the asymptotic of the number of
maps on a surface of genus $g$ arises from maps with schemes that have
the maximum number of vertices and edges. More precisely with
probability going to 1 as $n$ goes to infinity, a random map of genus
$g$ has a scheme with $4g-2$ vertices of degree 3.

The asymptotic behavior appears to have an interesting structure. As
shown in \cite{Chapuy} the more general construction of Bouttier et al
\cite{BDFG} in terms of labeled mobiles extends to higher genus in a
similar way as our original quadrangular construction. As a result one
obtains a more general class of labeled $g$-trees, say
$g$-mobiles. The decomposition of \cite{Chapuy} is more involved but
leads to the \emph{same dominant standard schemes}, with weights involving,
instead of $U$ and $T$, a finite number of algebraic series depending
on the set of allowed degrees.  In particular the same constant
\[
\tau_g=\sum_{\vec\ww\in\mathcal{W}^m_g}\prod_{i=1}^{4g-3}\frac1{d(i)}
\] 
drives the multiplicative constant in the asymptotic behavior, in
agreement with the pattern that was observed by \cite{Gao}, and with the
universality described in the physics literature \cite{DFGiZJ}.

\begin{acks}
This paper is based in part on early unpublished works of the last two
authors \cite{Ma:PhD,Sc:PhD} in their PhD theses under the direction
of Robert Cori.  We wish here to express our warm thanks to Robert
Cori for his interest and support for this project. We also thank
Gregory Miermont for communicating us \cite{Mi07} and for
interesting discussions.
\end{acks}
 
\bibliographystyle{alpha}
\bibliography{ChMaSc}
\label{sec:biblio}

\appendix 

\section{Elementary map surgery}\label{sec:surgery}
\subsection{Corners, edge addition and deletion} 

Let $v$ be a vertex of degree $k$ of a map $\mm$. The $k$ arcs with
origin $v$ in $\mm$ form a cycle around $v$ and the \emph{corner}
$c_{aa'}$ between two successive arcs $a$ and $a'$ in counterclockwise
order around $v$ is the sector of a small neighborhood of $v$ between
these arcs. The corner $c_{aa'}$ is said
to be \emph{incident} to $v$ and to the face which is on the left hand
side of $a$ (and on the right hand side of $a'$). A corner $c$ is
incident to a unique vertex $v$ and a unique face $f$, so we sometime
simply say that $c$ is a corner of $v$ or a corner of $f$.

Given a map $\mm$ and two (possibly identical) corners $c_1$ and $c_2$
incident to a same face $f$, there is a unique way to create a map
$\mm'$ by adding inside $f$ an edge $e$ that connects the vertices
incident to $c_1$ and $c_2$ at these corners: by this operation the
face $f$ of $\mm$ is naturally divided in $\mm'$ into two faces $f'$
and $f''$ separated by the edge $e$, and the corners $c_1$ and $c_2$
are respectively divided into $c'_1$ and $c''_1$ and $c'_2$ and
$c''_2$ accordingly. Observe that if $\mm$ is a map of genus $g$ with
$n$ edges, $v$ vertices and $f$ faces, then $\mm'$ is a map of genus
$g$ with $n+1$ edges, $v$ vertices and $f+1$ faces, in accordance with
Euler's formula.

Given a map $\mm$ and a subset $\ee$ of its edge set, we now consider
the problem of deleting $\ee$ from $\mm$. Let us denote by $\mm-\ee$
the graph drawing obtained by ignoring the edges of $\ee$ in the map
$\mm$. \emph{A priori}, $\mm-\ee$ can only be considered as a graph
drawing because the removal of some edges from $\mm$ can disconnect
the underlying graph or more generally it can create non simply
connected faces.  The following lemma gives sufficent conditions on
$\ee$ that garanty that $\mm-\ee$ is a map. From a topological point
of view, it is an elementary application of Van Kampen theorem.
\begin{lemma}\label{lem:surgery}
Let $f_1$ and $f_2$ be 2 adjacent faces of a map $\mm$ on $\mathcal
S_g$, and $e$ an edge separating $f_1$ and $f_2$. Then $\mm-\{e\}$ is a
 map on $\mathcal S_g$.

Let $v_0$ be a vertex of degree $d$ of a map $\mm$ on $\mathcal S_g$
such that $v_0$ is incident to $d$ different faces, and let $e_1,\ldots,e_d$
be the edges incident to $v_0$. Then $\mm-\{v_0,e_1,\ldots,e_d\}$ is a
map on $\mathcal S_g$.
\end{lemma}
\begin{proof}
For the first alinea, Van Kampen theorem can be applied directly: the
union of two simply connected regions with a simply connected
intersection is simply connected. The second alinea follows by
iteration.
\end{proof}


Implicitely in the previous discussion and in the rest of the text we
use the fact that the deletion or adjunction of edges and vertices to
a graph drawing that preserve the surface and the map structure
commute to the homeomorphisms of the surface, so that these operations
make sense up to homeomorphisms. As already said, the interested
reader can refer to \cite{MoTo} for a purely combinatorial treatment
of these questions.

\subsection{Duality and spanning trees}

Given a map $\mm$, the \emph{dual map} is obtained by exchanging the
role of vertices and faces and keeping incidence relations: more
precisely to construct the dual map $\mm^*$ of $\mm$, put a vertex in
each face of $\mm$ and connect these vertices by a dual edge $e^*$
across each edge of $\mm$. Given a map $\mm$ and a set of edges $\ee$
of $\mm$, let $\ee^*$ denote the subset $\{e^*\mid e\in\ee\}$ of the
set of edges of the dual map $\mm^*$. In terms of dual edges, the
first alinea of Lemma~\ref{lem:surgery} admits an immediate extension:
\begin{lemma}\label{lem:extendedsurgery}
Let $\ee$ be a set of edges of a map $\mm$ such that $\ee^*$ is a
forest (that is, the edges of $\ee^*$ do not form any cycle in
$\mm^*$). Then $\mm-\ee$ is a map on $\mathcal S_g$.
\end{lemma}
\begin{proof}
Let $e\in \ee$. We wish to apply Lemma~\ref{lem:surgery} to show that
$m- \{e\}$ is a map and conclude by induction since
$(\ee\setminus\{e\})^*$ is a forest with less edges. The only problem
could be that $e$ is incident twice to the same face $f$: but this
would mean, in dual terms, that $e^*$ is incident twice to the same
vertex, a contradiction with the assumption that the edges of $\ee^*$
form no cycle.
\end{proof}

We shall use the following immediate corollary of the previous lemma.
\begin{corollary}
If $\ee^*$ is a spanning tree of $\mm^*$ (that is, the edges of
$\ee^*$ do not form any cycle, and any two vertices of $\mm^*$ can be
joined by a path made of edges of $\ee^*$), then $\mm-\ee$ is a
$g$-tree.
\end{corollary}
\begin{proof}
By Lemma~\ref{lem:extendedsurgery}, $\mm-\ee$ is a map on
$\mathcal{S}_g$. Let $f_1$ and $f_2$ be two faces of $\mm-\ee$. Each
of them contains at least one vertex of $\mm^*$, and since $\ee^*$ is
a spanning tree these two vertices are connected by a path made of
edges of $\ee^*$. By definition of dual edges this path does not cross
the graph $\mm-\ee$: it stays in the same connected component of
$\mathcal{S}_g\setminus (\mm-\ee)$, so that $f_1=f_2$.
\end{proof}
Observe that in the planar case $g=0$, this is the classical duality result:
if $\ee^*$ spans $\mm^*$ then $\mm-\ee$ spans $\mm$.

\section{The correctness of the main bijection}\label{sec:proof}
\subsection{The mapping $\phi$ produces a well defined map}
Our first aim is to prove Lemma~\ref{lem:open}: the graph drawing
created by the opening $\phi$ is a map on $\mathcal S_g$ with one
face. This will essentially be done by applying
Lemma~\ref{lem:extendedsurgery}.

Let $(\qq,v_0)$ be a pointed bipartite quadrangulation of genus $g$
with $n$ faces, and let $\qq'=\qq\cup\phi(\qq,v_0)$: $\qq'$ is the
graph drawing obtained from $\qq$ by adding $n$ \emph{new} edges, one
in each face of $\qq$ as in the construction of $\phi(\qq,v_0)$. The
edges of $\qq$ remain in $\qq'$ and are called the \emph{old}
edges. Clearly, in view of the discussion of edge addition in
Section~\ref{sec:surgery}, $\qq'$ is a map of genus $g$: each face of
$\qq$ is divided into two faces in $\qq'$, either both of degree 3, or
one of degree 4 and one of degree 2. Then $\phi(\qq,v_0)$ is obtained
from $\qq'$ by deleting the old edges and the vertex $v_0$. We have to
check that deleting these old edges do not create a non simply
connected face.

Since $\qq'$ is a map, we can consider its dual $\qq''=(\qq')^*$. Let
us orient the edges of $\qq''$ that are duals of the old edges of
$\qq'$: an old edge of $\qq''$ whose dual edge has label $(i,i+1)$ is
oriented so that it has the largest label $i+1$ on its right hand side
(see Figure~\ref{fig:preuveduale}). The \emph{fake} edges of the proof
after Lemma~\ref{lem:open} are exactly these oriented dual edges.

\begin{lemma}\label{lem:dual}
There is a unique outgoing edge leaving each vertex of $\qq''$, so that
any cycle of oriented edges in $\qq''$ is in fact an oriented cycle.
The only oriented cycle in $\qq''$ consists of the duals of the edges
incident to $v_0$ in $\qq'$.
\end{lemma}
\begin{proof}
The first property follows immediately upon comparing our choice of
orientation with the three types of faces of $\qq'$, as represented on
Figure~\ref{fig:preuveduale}.

Now label each face of $\qq'$ by the minimum of the labels of its
incident vertices. The three cases represented in
Figure~\ref{fig:preuveduale} show that an oriented edge of $\qq''$
that crosses an edge of $\qq'$ with label $(i,i+1)$ originates in a
face with label $i$ and ends in a face with label $i$ or $i-1$: hence
an oriented cycle of edges of $\qq''$ must visit faces with non
increasing labels, which are thus all equal. A closer look at the
three cases furthermore shows that the only way to keep constant
labels consists in turning counterclockwise around a fixed vertex of
$\qq'$. An oriented cycle visiting faces with label $i$ in $\qq''$
must therefore enclose exactly one vertex. This vertex has then label
$i$ and all incident edges have label $(i,i+1)$. But, by definition of
the distance labeling, there is only one vertex not incident to a
vertex with a smaller label, the basepoint $v_0$.
\end{proof}

Lemma~\ref{lem:dual} implies that the oriented edges of $\qq''$ form
around $v_0$ an oriented cycle, to each vertex of which is attached a
tree of edges oriented toward the cycle (see
Figure~\ref{fig:preuveduale}).

Let us now consider the deletion of the old edges from $\qq'$.  Since
the possible labels around a face of $\qq'$ are $(i,i+1)$,
$(i,i+1,i+1)$ or $(i,i+1,i+2,i+1)$ for some $i\geq0$, no face of
$\qq'$ can be incident twice to $v_0$. Hence, according to the second
alinea of Lemma~\ref{lem:surgery}, removing the (old) edges incident
to $v_0$ yields a simply connected face. Keeping track of this
operation on the dual $\qq''$ amounts to the contraction of the
oriented cycle into a single vertex: the remaining oriented edges of
this modified dual then form a tree. According to
Lemma~\ref{lem:extendedsurgery}, removing these remaining old edges
from $\qq'$ yields again a simply connected face. This proves that
$\phi(\qq,v_0)$ is a well defined map on $S_g$, and that it has
exactly one face: $\phi(\qq,v_0)$ is a $g$-tree with $n$
edges. Finally, the labels are positive integers and the new edges
satisfy by construction the small variation condition so that
$\phi(\qq,v_0)$ is a well labelled $g$-tree.

\subsection{The closure of a well labeled $g$-tree}

We describe here a mapping $\psi'$, \emph{a priori} different from
$\psi$, but which we shall later prove to be the inverse of
$\phi$. While the description of $\psi$ relied on
Lemma~\ref{lem:luka}, we explicitely describe how to add edges to a
well labeled tree $\ttt$ to produce its image $\psi'(\ttt)$. In
particular the forthcoming discussion implicitely contains a
constructive proof of Lemma~\ref{lem:luka}, and shows that
$\psi'=\psi$.

Given a well labeled map, let $c$ be a corner with label $\ell$ of a
face $f$, and assume the minimal label around $f$ is
$\ell_0<\ell$. The \emph{predecessor} of $c$ is then defined as the
first corner with label $\ell-1$ encountered after $c$ when going
clockwise around the face $f$: this definition makes sense as
soon as the labels can decrease at most by one around any face, as is
the case for well labeled maps. In particular in a well labeled
$g$-tree, the predecessor is defined for any corner with label at
least 2, since there is a unique face with minimum label 1.

Let us now define the image $\psi'(\ttt)$ of a $g$-tree $\ttt$ with $n$
edges, and for this, assume that $\ttt$ has $k$ corners with label
$1$:
\begin{itemize}
\item[\emph{Phase 1.}]  Add a new vertex $v_0$ in the unique face of
 $\ttt$, and draw an edge between each corner of label $1$ of $\ttt$
 and $v_0$, so that the unique face of $\ttt$ is divided into $k$
 simply connected faces, each incident once to $v_0$. The result is a
 map $\ttt'$ of genus $g$ with $k$ faces, denoted $f_1,\ldots,f_k$ in
 clockwise order around $v_0$. Let also $\vec e_i$ denote the arc with
 origin $v_0$ that has $f_i$ on its right-hand side.
\end{itemize}
By construction, each face $f_i$ of $\ttt'$ is incident once to $v_0$,
it has one corner with label 0. All corners of $\ttt'$ thus have a
predecessor, except those incident to $v_0$.

\begin{itemize}
\item[\emph{Phase 2.}]  In each face $f_i$ of $\ttt'$, insert
 iteratively a edge between each corner and its predecessor, in the
 counterclockwise order of corners around $f_i$, starting with the
 first corner after $v_0$. The result is a map $\ttt''$ of genus $g$,
 containing the $n$ edges of $\ttt$ and $2n$ new edges.
\end{itemize}
As soon as an edge is added in $f_i$, this face is divided in the
current map $\mm$ into smaller faces, and some corners are divided
into smaller corners. However we will check that, due to the order of
insertions, the corner $c$ of $f_i$ is still a corner of $\mm$ when it
is considered, and that $\mm$ remains well labeled during all the
insersion process. It thus makes sense to ask for the insertion of an
edge between $c$ and its predecessor in $\mm$. Moreover the
predecessor of $c$ in $\mm$ is a corner which can be naturally
identified with the predecessor of $c$ in $\ttt'$, allowing for the
shorthand statement at the beginning of Step~2.

\begin{itemize}
\item[\emph{Phase 3.}] Remove the edges of the map $\ttt'$ from the map
$\ttt''$ to get a pointed graph drawing $\psi'(\ttt)=(\ttt''-\ttt',v_0)$.
\end{itemize}

\begin{lemma} The pointed graph drawing $\psi'(\ttt)$ is a pointed
bipartite quadrangulation on $S_g$ with $2n$ edges, and:
$$ \phi(\psi'(\ttt))=\ttt
$$
\end{lemma}

\begin{proof}
Let us first check that the construction is well defined.  Let $f_i$
be one of the $k$ faces of $\ttt'$, let $c_0,c_1,c_2,\ldots, c_q$ be
its corners in clockwise direction, starting with the corner of
$v_0$, and let $\ell_i$ be the label of $c_i$. At the $p$th step of
the insertion in $f_i$, the edges associated with corners
$c_1,\ldots,c_{p-1}$ have been inserted and the current map $\mm$
satisfy the following property, which is the invariant of
our construction:
\begin{itemize}
\item Let $f$ be the face on the right hand side of the arc $\vec e_i$ in
$\mm$: the corners of $f$ in counterclockwise order starting at $v_0$
are $c'_0,c'_1,c'_2,\ldots,c'_{\ell_{p-1}},c_p,c_{p+1},\ldots,c_q$, where
each $c'_i$ is a corner with label $i$ for $i=0,\ldots,{\ell_{p-1}}$. 
\end{itemize} 
In particular the corner $c_p$ of $f_i$ is a corner of $f$ in $\mm$
and its predecessor is the corner $c'_{\ell_{p}-1}$ (beware that
$c'_{\ell_p-1}$ is different from $c'_{\ell_{p-1}}$ unless
$\ell_{p-1}=\ell_p-1$).  The insertion of the edge $e$ between $c_p$ and
its predecessor divides $f$ into 2 faces:
\begin{itemize}
\item On the right hand side of $\vec e$ (oriented from $c_p$ to
$c'_{\ell_{p}-1}$), remains a face incident to the edge $e_i$ with corners
$c'_0,\ldots,c'_{\ell_p-2},c'',c''',c_{p+1}.\ldots,c_q$ where $c''$
and $c'''$ have respectively label $\ell_p-1$ and $\ell_p$. This ensures
the correctness of the invariant.
\item On the left hand side of $\vec e$, remains a face of degree 2,
3, or 4, depending on the value of $\ell_p-\ell_{p-1}$ (that is, on the
variation on the edge between corners $p-1$ and $p$ along
$f_i$). These 3 types of faces are illustrated by
Figure~\ref{fig:typesaretes}.
\end{itemize}
This concludes the proof that the construction is well defined.

In view of the 3 types of faces of $\ttt''$, it is then clear that the
edges of $\ttt$ with label $(i,i)$ lie in faces of $\psi'(\ttt)$ with
labels $i,i-1,i,i-1$, and that the edges of $\ttt$ with label
$(i,i+1)$ lie in faces of $\psi'(\ttt)$ with labels $i+1,i,i-1,i$ (see
Figure~\ref{fig:typesaretes}).  In particular each edge of $\ttt$ lies
in a face of $\psi'(\ttt)$ of degree 4. Hence $\psi'(\ttt)$ is a
quadrangulation, and it is bipartite because new edges always connect
vertices of different parity.

Upon comparing the two cases in Figure~\ref{fig:typesaretes} with the
two cases for the creation of edges by $\phi$ in
Figure~\ref{fig:rules}, it clearly appears that in any face of
$\psi'(\ttt)$ the construction rules of $\phi(\psi'(\ttt))$ recreates the
original edge of $\ttt$, so that $\phi(\psi'(\ttt))=\ttt$.
\end{proof}

\subsection{Proof of the main theorem}

To complete the proof of Theorem~\ref{thm:bijection}, it is sufficient
to show that for every pointed bipartite quadrangulation $(\qq,v_0)$ and
associated $g$-tree $\ttt=\phi(\qq,v_0)$, one has $\psi'(\ttt)=(\qq,v_0)$.  The
key step is the following lemma:
\begin{lemma}
\label{lemma:last}
Let $(\qq,v_0)$ be a pointed bipartite quadrangulation, let
$\ttt=\phi(\qq,v_0)$ be the associated $g$-tree. Let $\qq'=\qq\cup\ttt$,
that is $\qq'$ is the map in which the edges of $\ttt$ have been added
to $\qq$, and let $\ttt''=\ttt\cup\psi'(\ttt)$, that is the map in
which the edges of $\psi'(\ttt)$ have been added to $\ttt$. 
Both $\qq'$ and $\ttt''$ contain the $g$-tree $\ttt$.

If there is an edge $e$ of $\qq$ that connects two corners of $\ttt$
in $\qq'$ (resp. a corner of $\ttt$ and $v_0$), then there is in $\ttt''$ an 
edge between the same corners of $\ttt$ (resp. between the same corner of 
$\ttt$ and $v_0$).
\end{lemma}

\begin{figure}
\centerline{\includegraphics[scale=1]{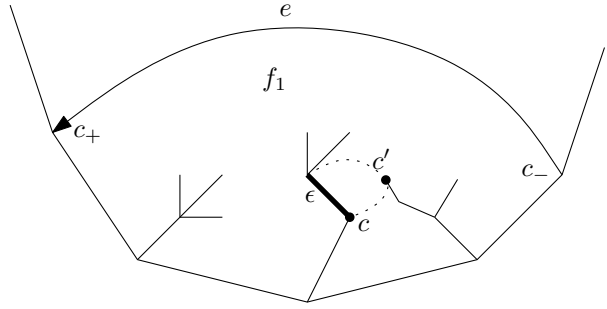}}
\caption{The proof of Lemma~\ref{lemma:last}}\label{fig:lastfig}
\end{figure}

\begin{proof}
If $e$ is adjacent to $v_0$, then it is immediately replaced in
$\psi'(\ttt)$ during the first phase of the construction, and the conclusion of the
lemma holds.

Let us thus consider an edge $e$ of $\qq$ not adjacent to
$v_0$, and let us show that it is recovered in $\ttt''$.
Let $\ttt'$  be the map obtained by adding to $\ttt$ the vertex $v_0$
and all the edges adjacent to $v_0$ in $\qq$: $\ttt'$
is a submap of $\qq'$, and each face of $\ttt'$ is adjacent to exactly one 
corner with label $0$.
Since $\ttt'$ spans all the vertices of $\qq'$, $e$ must lie in
some face $f$ of $\ttt'$, cutting it into two new faces. Call $f_1$
the one of these two faces which does not contains the corner labeled
$0$, and orient $e=(c_-,c_+)$ in such a way that it has $f_1$ on its
left. Here, $c_-$ and $c_+$ are two corners of $\ttt$, but up
to a slight abuse of notations, can be considered as corners of $f_1$
as well.

We claim that $c_+$ is the unique corner of minimal label in
$f_1$. Indeed, assume that the minimal label in $f_1$ is reached
at some corner $c \neq c_+$, and let $\ell_c$ be its label.
Let $\epsilon$ be the edge of $\ttt$ leaving $c$ in clockwise direction 
around $f_1$. By the rules of construction of $\ttt$, $\epsilon$ has been 
obtained from
a face of $\qq$ that contains a vertex $c'$ of label $\ell_c-1$ on the
right of $\epsilon$. Since $\qq'$ is a well defined map,
with no vertex strictly inside $f_1$, 
the vertex $c'$ must lie on the border of $f_1$ (see Figure~\ref{fig:lastfig}).
This contradicts the minimality of the label $\ell_c$ in $f_1$.

Hence $c_+$ is the unique corner of minimal label in $f_1$. 
Since $e$ is an edge of $\qq$ this implies that $\ell_{c_-}=\ell_{c_+}+1$, and
that $c_+$ is the predecessor of $c_+$ in $\ttt$: $c_+$ and $c_-$
are linked by an edge in $\ttt''$.
\end{proof}

From the lemma, the set of edges of $\qq$ is naturally included in the set of 
edges of $\psi'(\ttt)$.
Now, by construction, the number of edges in $\psi'(\ttt)$ is equal to 
the number of edges in $\qq$. Hence $\psi'(\ttt)$ and $\qq$ 
have exactly the same set of edges, and thus are equal. 
This completes the proof of Theorem~\ref{thm:bijection}.

\end{document}